\documentclass[12pt,leqno]{article}
\usepackage{amsmath}
\usepackage{amssymb}
\usepackage{theorem}

\numberwithin{equation}{section}

\def\picture #1 by #2 (#3){
  \vsquare to #2{
    \hrule width #1 height 0pt depth 0pt
    \hfill
    \special{picture #3}}}

\def\scaledpicture #1 by #2 (#3 scaled #4){{
  \dimen0=#1 \dimen1=#2
  \divide\dimen0 by 1000 \multiply\dimen0 by #4
  \divide\dimen1 by 1000 \multiply\dimen1 by #4
  \picture \dimen0 by \dimen1 (#3 scaled #4)}}

\theoremstyle{plain}
\newtheorem{thm}{Theorem}

\newtheorem{lem}[thm]{Lemma}

\def\D{\Delta}

\def\Dr{\D_\rho}

\def\supp{\mathrm{supp }}
\renewcommand{\span}{\mathrm{span}}

\def\dist{\mathrm{dist}}
\def\R{\mathbb{R}}
\def\cS{\mathcal{S}}

\newcommand{\po}{\partial\Omega}

\def\C{\mathbb{C}}

\def\dbc{\overline\partial \chi_0}

\def\Q{{\mathcal Q}}

\def\U{{\mathcal U}}

\def\uo{u^{(1)}}
\def\ut{u^{(2)}}

\def\CD#1{\mathcal C\mathcal D_{q_{#1}}}
\def\F#1{\mathcal F_{q_{#1}}}
\def\u#1{u^{(#1)}}
\def\v#1{v^{(#1)}}

\title{Local uniqueness for the Dirichlet-to-Neumann map via the two-plane
transform}
\author{Allan Greenleaf\thanks{Partially supported by NSF grant DMS-9877101.}\\
Department of Mathematics\\
University of Rochester\\
Rochester, NY 14627
\and
Gunther Uhlmann\thanks{Partially supported by NSF
grant DMS-9705792 and the Royal Research Fund at the University of Washington}\\
Department of Mathematics\\
University of Washington\\
Seattle, WA 98195}
\date{Revised - November 8, 2000}

\begin{document}

\maketitle

\begin{abstract} We consider the Cauchy data
associated
to the Schr\"odinger equation with a potential on a bounded domain $\Omega
\subset \R^n, n\ge 3$.
We show that the integral  of the potential over a two-plane $\Pi$ is
determined by the Cauchy data of certain exponentially
growing solutions on any
 open subset $\mathcal
U\subset\partial\Omega$ which contains
$\Pi\cap\partial\Omega$.
\end{abstract}

\setcounter{section}{-1}
\section{Introduction}

For $\Omega$ a bounded domain in $\R^n$ with Lipschitz boundary,
$\partial\Omega$, and real-valued $q(x)\in L^\infty(\Omega)$ , let
\begin{equation}\label{eqn0.1}
\Lambda_q:H^{\frac{1}{2}}(\partial\Omega)\to
H^{-\frac{1}{2}}(\partial\Omega)
\end{equation}
be the Dirichlet-to-Neumann map associated with the operator $\Delta+q$ on
$\Omega$, which is defined if $\lambda=0$ is not a Dirichlet
eigenvalue
for $\Delta +q$ on $\Omega$.
More generally, one may consider the set of Cauchy data of solutions
of $(\D+q(x))v=0$, which is defined even if $\lambda=0$ is a Dirichlet
eigenvalue. Set
\begin{equation}
\CD{ }=\Bigl\{(v|_{\po},\frac{\partial v}{\partial n}|_{\po})
\in H^{\frac12}(\po)\times H^{-\frac12}(\po):
v\in H^1(\Omega), (\D+q)v=0\Bigr\},
\end{equation}
which is a subspace of $H^{\frac12}\times H^{-\frac12}$ ;
if $\Lambda_q$ is defined,
then $\CD{ }$ is simply the graph of $\Lambda_q$.

This paper is concerned with
the problem of obtaining partial knowledge of $q(x)$ from partial knowledge
of  $\CD{ }$, namely its restriction to certain
``small'' open subsets of
the boundary.   The approach taken here is
to  use concentrated,  exponentially growing, approximate solutions to
relate
$\CD{ }$ on an open set $\mathcal{U}\subset\po$ to the two-plane transform
of
the potential
$q(x)$ on two-planes whose intersections with $\po$ are contained in
$\mathcal{U}$.

Let $M_{2,n}$ denote the
$(3n-6)$-dimensional Grassmannian of all affine two-planes $\Pi\subset\R^n$,
and \begin{equation}\label{eqn0.3}
R_{2,n}f(\Pi)=\int_\Pi f(y)d\lambda_\Pi(y), f\in L^2(\R^n),
\end{equation}
denote the two-plane transform on $\R^n$ \cite{H65, H80}. Here, $d\lambda_\Pi$
is two-dimensional Lebesgue measure on $\Pi\in M_{2,n}$, which can be defined
by
\begin{equation}
<f,d\lambda_\Pi>=\lim_{\epsilon\rightarrow
0}\frac{1}{|B^{n-2}(0;\epsilon)|}\int_{\{dist(x,\Pi)<\epsilon\}}f(x) dx.
\end{equation}
(Note
that for $n=3$, $R_{2,3}$ is just the usual Radon transform on $\R^3$.)
We will also need the variant of $d\lambda_\Pi$ defined relative to
$\Omega$:
\begin{equation}
<f,d\lambda_\Pi^\Omega>=\lim_{\epsilon\rightarrow
0}\frac{1}{|B^{n-2}(0;\epsilon)|}\int_{\Omega\cap\{dist(x,\Pi)<\epsilon\}}f(
x)
dx,
\end{equation}
which gives rise to a two-plane transform relative to $\Omega$,
\begin{equation}\label{eqn0.4}
R_{2,n}^\Omega f(\Pi)=\int_\Pi f(x)d\lambda_\Pi^\Omega(x).
\end{equation}
Note that if $\po$ is $C^1$ and $\Pi\cap\po$ transversally, then
$<d\lambda_\Pi^\Omega,f>=\linebreak<d\lambda_\Pi,f\cdot\chi_\Omega>$ and
$R_{2,n}^\Omega f(\Pi)=R_{2,n}(f\cdot\chi_\Omega)(\Pi)$.

For each choice of an orthonormal basis for $\Pi_0$, the translate of $\Pi$
passing
through the origin, as well as other arbitrary choices made below, we will
construct a
family, $\F{ } =\{v_z(x): z\in\C, |z|\ge C\}$, of exponentially growing
solutions of
$(\D+q(x))v=0$, concentrated near $\Pi$. Using these families, we formulate
\medskip

{\bf Definition}
(i) If $\U\subset\po$ is open,  $\CD{1}$ and $\CD{2}$ are
{\it equal on $\U$ relative to $\mathcal F$ at $z\in\C$} if the solutions in
$\F{1}$ and
$\F{2}$ corresponding to opposite exponential growths, $v^{(1)}_z$ and
$v^{(2)}_{-z}$,
have the same Cauchy data on $\U$:
\begin{equation}\nonumber
(v^{(1)}_z|_{\U},\frac{\partial v^{(1)}_z}{\partial n}|_{\U})=
(v^{(2)}_{-z}|_{\U},\frac{\partial v^{(2)}_{-z}}{\partial n}|_{\U}).
\end{equation}
(ii) $\CD{1}$ and $\CD{2}$ are {\it equal on $\U$ for a sequence of
exponentially
growing solutions} if $\CD1$ and $\CD2$ are equal on $\mathcal U$ relative
to
$\mathcal F$ at $z=z_j$ for
some sequence $\{z_j\}^\infty_1\subset\C$ with $|z_j|\rightarrow\infty$.

\medskip

We may now state the main result proved here.
For each $\Pi\in
M_{2,n}$, let  $\gamma_\Pi=\Pi\cap\po\subset\po$, and let $H^{s}(\Omega)$
denote the standard Sobolev space of distributions with $s$
derivatives in $L^2(\Omega)$.
\bigskip

\begin{thm}\label{thm1}
Let $n\ge 3$. Assume $\po$ is Lipschitz and potentials $q_1(x)$ and $q_2(x)$
are in  $H^{s}(\Omega)$, for  some $s>\frac{n}2$.
Let $\Pi\in M_{2,n}$ and $\F{1}$ and $\F{2}$ be families of exponentially
growing solutions associated to $q_1$ and $q_2$. If, for some fixed
neighborhood
$\U_\Pi$ of $\gamma_\Pi$ in $\po$, $\CD{1}$ and $\CD{2}$ are equal on
$\U_\Pi$ for a sequence of exponentially growing solutions,
then
\begin{equation}\label{eqn0.5}
R_{2,n}^\Omega(q_1-q_2)(\Pi)=0,
\end{equation}
i.e., $\int q_1(y)d\lambda_\Pi^\Omega (y) = \int
q_2(y)d\lambda_\Pi^\Omega (y)$.
\end{thm}

\bigskip If $\CD{1}$ and $\CD{2}$ equal on all of $\po$
relative to ${\mathcal F}$, then
 this implies that
$R_{2,n}((q_1-q_2)\chi_\Omega)(\Pi)=0$,
$\forall\;\Pi\in M_{2,n}$, which by the uniqueness theorem for $R_{2,n}$
yields that $q_1-q_2\equiv 0$ on $\Omega$, providing a variant of the
global uniqueness theorem for the Dirichlet--to--Neumann
map \cite{SU87a}. (We note that our technique is limited to three
or more
dimensions and says nothing in the case $n=2$ \cite{N96}.)  However,
one is also able
to obtain {\it local} uniqueness results by replacing the uniqueness theorem
for the two-plane transform with Helgason's support theorem [H80,
Cor.\ 2.8]: if $C\subset\R^n$ is a closed, convex set and $f(x)$ a
function\footnote{The support and uniqueness theorems are usually stated
under the assumption that $f(x)$ is continuous, of rapid decay in the case
of the support theorem, but the proofs in [H80] are easily seen to extend
to the case where $f(x)=q(x)\chi_\Omega(x)$ with $\Omega\subset\R^n$
bounded, $q\in C(\overline\Omega)$.} such that $R_{2,n}f(\Pi)=0$ for all
$\Pi$ disjoint from $C$, then $\supp(f)\subset C$. We then immediately
obtain the following two results.

\begin{thm}\label{thm2}
Suppose $\po$ and potentials $q_1,q_2$
are as in Thm. 1., and $C\subset\Omega$ is a closed, convex set.
 If, for all $\Pi\in M_{2,n}$ such that $\Pi\cap C=\phi$,
there is some neighborhood $\U_\Pi$ of
$\gamma_\Pi$ on which  $\CD{1}$ and $\CD{2}$ are
equal for some sequence of exponentially growing solutions,
then $\supp(q_1-q_2)\subseteq C$, i.e., $q_1=q_2$ on $\Omega\backslash
C$. \end{thm}

\begin{thm}\label{thm3}
Suppose $\po$ is $ C^{2}$ and strictly convex, and potentials $q_1,q_2$
are as in Thm. 1. If, for some $r>0$,
$\CD{1}$ and $\CD{2}$ are
equal on $B$ for some sequence of exponentially growing solutions
 for all surface balls
$B=B^n(x_0;r)\cap\po\subset\po$, then \[
\dist(\supp(q_1-q_2),\po)\ge Cr^2,
\]
i.e., $q_1=q_2$ on the tubular neighborhood $\{x\in\overline\Omega:\dist(x,
\po)\le Cr^2\}$ of $\po$ in $\overline{\Omega}$.
\end{thm}

\noindent {\bf Remark}

\par The conclusions of Thms. 2 and 3 can be strengthened by combining them
with a result in Isakov [Is]. Namely, if either $C\subset\subset\Omega$
in Thm. 2, or the assumption of Thm. 3 holds for some  $r>0$,
we can conclude from Thm. 2 or 3 that
$supp(q_1-q_2)\subset\subset\Omega$.
 By Ex. 5.7.4 in [Is], based on a technique of Kohn and Vogelius[KV85],
this, together with the condition that
$\Lambda_{q_1}=\Lambda_{q_2}$ on some
open set $\U\subset\po$,  implies that $q_1\equiv q_2$ everywhere
on $\Omega$. We are indebted to Adrian Nachman for pointing this out to us.

The authors would like to thank Alexander Bukhgeim and Masaru Ikehata for
pointing
out  errors in an earlier version of this paper.

\section{Approximate solutions}

To prove Thm.\ \ref{thm1}, we first construct  exponentially growing
{\it approximate} solutions for $(\D+q)v=0$. As considered
in \cite{C, SU86, SU87a}, let
\[
\Q = \{\rho\in\C^n:\rho\cdot\rho=0\}
\]
be the (complex) characteristic variety of $\Delta$. Each $\rho\in\Q$ can be
written as
$\rho=|\rho|\frac{\rho}{|\rho|}=\frac{1}{\sqrt{2}}|\rho|(\omega_R+i\omega_I)
\in\R\cdot
\left( S^{n-1}+i S^{n-1}\right)$, with
$\omega_R\cdot\omega_I=0$. For $\rho\in\Q$, let $\Delta_\rho=\Delta
+2\rho\cdot\nabla$. Then
\begin{equation}\label{eqn1.1}
\Delta_\rho +q(x)=e^{-\rho\cdot x}(\Delta+q(x))e^{\rho\cdot x},
\end{equation}
so that, with $v(x)=e^{\rho\cdot x}u(x)$,
\begin{equation}\label{eqn1.2}
(\Delta_\rho+q(x))u(x)=w(x)\Leftrightarrow (\Delta+q(x))v(x)=e^{\rho\cdot
x}w(x)
\end{equation}
and, in particular, $(\Delta_\rho+q(x))u(x)=0\Leftrightarrow
(\Delta+q(x))v(x)=0$.

Now, given a potential $q(x)$ and a two-plane $\Pi\in M_{2,n}$, we will
construct an approximate solution $u_{app}$ to $(\Dr+q)u=0$, supported
near
$\Pi$:

\begin{thm}\label{thm4}
Let $\Omega$ be Lipschitz and $q(x)\in H^{s}(\Omega)$ for some
$s>\frac{n}2$. Then,
 for any $0<\beta<\frac14$
fixed, the following holds: $\exists\;\epsilon>0$ such that, for any
$\rho=\frac{1}{\sqrt{2}}|\rho|(\omega_R+i\omega_I)\in \Q$
 and any two-plane $\Pi$ parallel to
$\Pi_0=\span\{\omega_R,\omega_I\}$, we can find an approximate solution
$u_{app}=u_{app}(x,\rho,\Pi)$ to $(\Delta_\rho+q(x))u=0$ satisfying
\begin{gather*}
\|u_{app}\|_{L^2(\R^n)}\le C,\quad
\|u_{app}\|_{L^2(\Omega)}\simeq[\lambda_\Pi^\Omega(\Pi\cap\Omega)]^{\frac{1}
{2}}\hbox{ as } |\rho|\rightarrow\infty
\tag{1.5} \label{eqn1.5} \\
\supp(u_{app})\subset\biggl\{x\in\R^n:\dist(x,\Pi)\le
\frac{2}{|\rho|^\beta}\biggr\}
\tag{1.6} \label{eqn1.6}
\intertext{and}
\|(\Delta_\rho+q)u_{app}\|_{L^2(\Bbb R^n)} \le
\frac{C_\epsilon}{|\rho|^\epsilon}.
\tag{1.7}
\label{eqn1.7}
\end{gather*}
Furthermore, for any two such solutions, $u_{app}^{(1)}, u_{app}^{(2)}$,
associated with possibly different potentials $q_1(x), q_2(x)$ and with
$\rho_1\in\Q,
\rho_2=e^{i\theta}\rho_1\hbox{ or
}\rho_2=e^{i\theta}\overline{\rho_1}\in\Q$,
\setcounter{equation}{7}
\begin{equation}\label{eqn1.8}
u_{app}^{(1)}(\cdot, \rho_1,\Pi) u_{app}^{(2)}(\cdot, \rho_2,
\Pi) \to
d\lambda_\Pi^\Omega\hbox{ weakly as }|\rho_1|\to\infty.
\end{equation}
\end{thm}
\bigskip
In fact, as will be seen below, $u_{app}=u_0+u_1$ with $u_0$ depending only
on
$\Pi$ and $|\rho|$ and satisfying (1.5).
\medskip

Now, we may apply the results of \cite{SU86,SU87a} (see
also \cite{H96}) to find a solution $u_2$ of
$$(\Dr+q)u_2=-(\Dr+q)u_{app}\in L^2_{comp}(\R^n),\notag$$
uniformly in $H^1_t$ and with a gain of $|\rho|^{-1}$ in $L^2_t$, as
long as $|\rho|\ge C$ with $C$ depending only on $\|q\|_\infty$
and $ diam(\Omega)$.
Here, $H^s_t$ and $L^2_t$  the
weighted versions
of these spaces, as in [SU87a], for some fixed $-1<t<0$.
By these results and (1.7),
$$\|u_2\|_{H^1_t(\R^n)}\le c
\|(\Dr+q)u_{app}\|_{L^2_{t+1}(\R^n)}\le c |\rho|^{-\epsilon},
\quad \|u_2\|_{L^2_t}\le C|\rho|^{-1-\epsilon}.\notag$$
(The statements in [SU86,SU87a] are for $q\in C^\infty$, but the proofs are
easily seen to hold if $q\in H^s(\Omega)$ with  $s>\frac{n}2$. Also,
the weights will be irrelevant since we will be working on $\Omega$.)
Thus, $u=u_{app}+u_2=u_0+u_1+u_2$ is an exact solution of $(\Dr+q)u=0$ on
$\R^n$,
satisfying
$$\|u-u_0\|_{L^2}\le c|\rho|^{-\epsilon}\hbox{ and } \|u_2\|_{H^s}\le
|\rho|^{s-1-\epsilon},\forall 0\le s\le 1.\notag$$

Finally,
$$\F{ }=\bigl\{v_z: |z|\ge C\bigr\}=
\bigl\{e^{\rho\cdot x}u(x,\Pi,\rho): \rho=
Re(z)\omega_R+i Im(z)\omega_I, |z|\ge C\bigr\}\notag$$
is the associated family of exponentially growing solutions used in the
statements of
the theorems.
\bigskip
To prove Thm. 1, we assume that $q_1,q_2$ and $\Pi\in M_{2,n}$,
$\U_\Pi\subset\po$ are as in its statement. We will make use of a
variant of Alessandrini's identity \cite{A}. For $j=1,2$, let
$v^{(j)}_{\rho_j}$
be the exact solution to $(\D+q_j)v=0$ constructed above, so that
$v^{(j)}_{\rho_j}(x)=e^{\rho_j\cdot x} u^{(j)}(x,\Pi,\rho_j)$, with
$u^{(j)}=u^{(j)}_{app}+u^{(j)}_2$. Taking
$\rho_1=\rho, \rho_2=-\rho$, consider the quantity
$$I=\int_{\po} \frac{\partial\v{1}_\rho}{\partial n}\cdot\v{2}_{-\rho}
-\v{1}_{\rho}\cdot\frac{\partial\v{2}_{-\rho}}{\partial n} d\sigma.\notag$$
Under the assumption that $\v{1}_\rho$ and $\v{2}_{-\rho}$ have the same
Cauchy data on
$\U_\Pi$, $I$ is equal to the integral of the same expression over
$\po\backslash\U_\Pi$. Observing that
$$\frac{\partial\v{1}_\rho}{\partial n}=e^{\rho\cdot
x}(\frac{\partial}{\partial n}
+(\rho\cdot n(x)))\u{1}\hbox{ and }
\frac{\partial\v{2}_{-\rho}}{\partial n}=e^{-\rho\cdot
x}(\frac{\partial}{\partial n}
-(\rho\cdot n(x)))\u{2},\notag$$
we see that the exponentials cancel and the integrand of $I$ is
$$=\frac{\partial \u{1}}{\partial n}\cdot\u{2}-\u{1}\cdot
\frac{\partial \u{2}}{\partial n} +2(\rho\cdot n(x))\u{1}\u{2}.\notag$$
Since (1.6) implies that $\supp(\u{j}_{app}|{\po}),
\supp(\frac{\partial\u{j}_{app}}
{\partial n}|_{\po})
\subset\U_\Pi$ for $|\rho|$ sufficiently large, we have that
$$I=\int_{\po\backslash\U_\Pi}
\frac{\partial \u{1}_2}{\partial n}\cdot\u{2}_2-\u{1}_2\cdot
\frac{\partial \u{2}_2}{\partial n} +2(\rho\cdot n(x))\u{1}_2\u{2}_2
d\sigma.
\notag$$
We estimate
\begin{eqnarray*}
|\int_{\po\backslash\U_\Pi} \frac{\partial \u{1}_2}{\partial n}\cdot\u{2}_2
d\sigma | & \le & \|  \frac{\partial \u{1}_2}{\partial n} \|_{H^{-\frac12}(
\po)}\cdot \| \u{2}_2 \|_{H^{\frac12}(\po)} \\
 & \le & \| \u{1}_2 \|_{H^{\frac12}(\po)} \cdot \| \u{2}_2
\|_{H^{\frac12}(\po)}\\
 & \le & C \| \u{1}_2 \|_{H^{1}(\Omega)} \cdot \| \u{2}_2
\|_{H^{1}(\Omega)}\hbox{ by Sobolev restriction }\\
 & \le & C \| \u{1}_2 \|_{H^{1}_t(\R^n)} \cdot \| \u{2}_2
\|_{H^{1}_t(\R^n)}\hbox{ since }\Omega\hbox{ compact }\\
 & \le & C|\rho|^{-2\epsilon}\rightarrow 0\hbox{ as }
|\rho|\rightarrow\infty
\end{eqnarray*}
and similarly for the second term. Now note that $|\rho\cdot n(x)|\le
c|\rho|
$ since $\po$ is Lipschitz, and
$$\|\u{j}_2\|_{L^2(\po)}\le \|\u{j}_2\|_{H^{\sigma}(\po)}
\le c_{\sigma} \|\u{j}_2\|_{H^{\sigma+\frac12}(\Omega)}
\le c_{\sigma}'|\rho|^{\sigma-\frac12-\epsilon}\notag$$
for any $\sigma>0$, and thus the third term is dominated by
$(c_{\sigma}')^2|\rho|\cdot |\rho|^{2\sigma-1-2\epsilon}\rightarrow 0$
as $|\rho|\rightarrow 0$ if we choose $0<\sigma<\epsilon$.

On the other hand,

\begin{eqnarray*}
I & =  & \int_{\po} \frac{\partial\v{1}}{\partial
n}\cdot{\v{2}}-\v{1}\cdot{\frac{\partial\v{2}}{\partial n}}
d\sigma \\
& = & \int_\Omega \D(\v{1})\cdot{\v{2}}-\v{1}\cdot{\D(\v{2})} dx
\hbox{
by Green's Thm. } \\
& = & \int_\Omega (
-q_1\v{1})\cdot{\v{2}}-\v{1}\cdot{(-q_2\v{2})}
dx \\
& = & \int_\Omega
(q_2-q_1)\v{1}{\v{2}}dx = \int_\Omega
(q_2-q_1)\u{1}{\u{2}}dx
\end{eqnarray*}
since the exponentials cancel. As $\u{1}\cdot\u{2}=(\uo_{app}+\uo_2)\cdot
(\ut_{app}+\ut_2)$ and the leading term $\uo_{app}\ut_{app}\rightarrow
d\lambda^{\Omega}_{\Pi}$ weakly as $|\rho|\rightarrow\infty$ by (1.8),
while the remaining terms $\rightarrow 0$ since
$\|\u{j}_{app}\|_{L^2(\Omega)}
\le C$ by (1.5) and $\|\u{j}_2\|_{L^2(\Omega)}\le c|\rho|^{-1-\epsilon}$,
we conclude that
$I\rightarrow
R^\Omega_{2,n}\bigl(q_2-q_1\bigr)(\Pi)\hbox{ as }
|\rho|\rightarrow\infty$,
finishing the proof of Thm. 1.
\bigskip

Now, to start the proof of Thm.\ \ref{thm4} we may use the rotation
invariance
of
$\Delta$ and the
invariance of $\Q$ under $S^1=\{e^{i\theta}\}$,  and note that it suffices
to treat the
case\footnote{Of course, the length of this element of $\Q$ is
$\sqrt{2}|\rho|$, but this is irrelevant for the proofs, and 
denoting the length of $|\rho|(\vec{e_1}+i\vec{e_2})$ 
by $|\rho|$ is notationally
convenient.} $\rho=|\rho|(\vec{e_1}+i\vec{e_2})$, where $\{\vec{e_1},\ldots,
\vec{e_n}\}$ is the standard orthonormal basis for $\R^n$. Write $x\in\R^n$
as $x=(x', x'')\in\R^2\times\R^{n-2}$ and similarly $\xi=(\xi',\xi'')$.

If $\Pi\in M_{2,n}$ is parallel to $\span\{\omega_R,
\omega_I\}=\span\{\vec{e_1}, \vec{e_2}\}=\R^2\times\{0\}$, then
$\Pi=\span\{\vec{e_1}, \vec{e_2}\}+(0, x''_0)$ for some $x''_0\in\R^{n-2}$.
Given $|\rho|>1$ and $x''_0\in\R^{n-2}$, we will define an approximate
solution $u(x, \rho, \Pi)$ to $(\Delta_\rho+q(x))u=0$ on $\R^n$, of the
form $u(x, \rho, \Pi)=u_0(x, \rho, \Pi)+u_1(x, \rho, \Pi)$.

For notational
convenience, we will usually suppress the dependence on $\rho$ and $\Pi$ and
simply write $u(x)=u_0(x)+u_1(x)$. We will use various cutoff functions
$\chi_j$;
for $j$ even or odd, $\chi_j$ will always denote a
function of $x'$ or $x''$, respectively.
Also, $B^m(a;r)$ and $S^{m-1}(a;r)$ will denote the closed ball
and sphere of radius $r$ centered
at a point $a\in\R^m$.

To define $u_0$, first fix $\chi_0\in C^\infty_0(\R^2)$ with $\chi_0\equiv
1$
on $B^2(0; R)$ for any $R>sup\{|x'|:(x', x'')\in\Omega$
for some $x''\in\R^{n-2}\}$; let $C_0=\|\chi_0\|_{L^2(\R^2)}$.
Secondly, let $\psi_1\in C^\infty_0(\R^{n-2})$ be radial, nonnegative,
supported in the unit ball, and satisfy
\begin{equation}\notag
\int_{\R^{n-2}}(\psi_1(x''))^2dx''=1.
\end{equation}

Now, for $\beta>0$ to be fixed later, we let $\delta$ be the small parameter
$\delta=|\rho|^{-\beta}$ and define
\[
\chi_1(x'')=\delta^{-\frac{n-2}{2}}\psi_1\left(\frac{x'-x''_0}{\delta}\right
),
\]
so that

\begin{equation}\label{eqn1.11}
\|\chi_1\|_{L^2(\R^{n-2})} =
\|\psi_1\|_{L^2(\R^{n-2})}=1, \;\forall\;\delta>0.
\end{equation}
Set
$u_0(x)=u_0(x', x'')=\chi_0(x')\chi_1(x'')$; then
$u_0$ is real, $\|u_0\|_{L^2(\R^n)}=C_0$
and $\|u_0\|_{L^2(\Omega)}\to [\lambda_\Pi(\Pi\cap\Omega)]^{\frac{1}{2}}$ as
$\delta\to 0^+$, i.e., as $|\rho|\to\infty$. Note
also that $||u_0||_{H^{1}}\le c
\delta^{-1} = c|\rho|^{\beta}$, so that
$||u_0||_{H^{s}}\le c |\rho|^{s\beta}$ for $0\le s\le 1$. Since
$\Delta_\rho=\Delta+2\rho\cdot\nabla=\Delta+2|\rho|
(\vec{e_1}+i\vec{e_2})\cdot\nabla=\Delta+4|\rho|\bar\partial_{x'}$ and
$\rho\bot\R^{n-2}$,
\begin{eqnarray*} (\Delta_\rho+q(x))u_0 & = &
(\Delta\chi_0)\cdot\chi_1+2(\nabla\chi_0)
\cdot(\nabla \chi_1)+\chi_0(\Delta\chi_1) \\
&& \quad\quad
+2(\rho\cdot\nabla)(\chi_0)\chi_1+2\chi_0(\rho\cdot\nabla)(\chi_1) +
q\chi_0\chi_1 \\
& = & \chi_0(x')(\Delta_{x''}+q)(\chi_1)(x'')\hbox{ on
}B^2(0;R)\times\R^{n-2},
\end{eqnarray*}
the first and fourth terms after the first equality vanishing because
$(\rho\cdot\nabla)(\chi_0)=2\overline\partial \chi_0\equiv 0$ on $B^2(0;R)$,
and the second and fifth equalling
zero because $\nabla \chi_1\bot\R^2$.

To define the second term in the approximate solution, $u_1(x)$, we make use
of a truncated form of the Faddeev Green function, $G_\rho$,  and an
associated
projection operator. The operator $\Delta_\rho$ has, for $\rho\in \Q$,
(full) symbol
\begin{equation}
\sigma(\xi)=-[(|\xi|^2-2|\rho|\omega_I\cdot\xi)+i2|\rho|(\omega_R\cdot\xi)],
\end{equation}
and so for $\frac{\rho}{|\rho|}=e_1+ie_2$, we have
\begin{equation}\notag
\sigma(\xi)=-[(|\xi-|\rho|\vec{e_2}|^2-|\rho|^2)+i(2|\rho|\xi_1)],
\end{equation}
which has (full) characteristic variety
\begin{eqnarray}\label{eqn1.12}
\Sigma_\rho & = & \{\xi\in\R^n:\xi_1=0, |\xi-|\rho|e_2|=|\rho|\} \\
& = & \{0\}\times S^{n-2}((|\rho|,0,\ldots,0);|\rho|)
\subset\R_{\xi_1}\times\R^{n-1}_{\xi_2,\xi''}.\nonumber
\end{eqnarray}
The Faddeev Green function  is then defined by
$G_\rho=(-\sigma(\xi)^{-1})^{\vee}\in\cS'(\R^n)$. We now introduce,
for an $\epsilon_0>0$ to be fixed later, a tubular
neighborhood of $\Sigma_\rho$,
\begin{equation}\label{eqn1.13}
T_\rho = \{\xi:\dist(\xi,\Sigma_\rho)<|\rho|^{-\frac{1}{2}-\epsilon_0}\},
\end{equation}
as well as its complement, $T^C_\rho$, and let $\chi_{T_\rho}$,
$\chi_{T_\rho^C}$
be their characteristic functions. Define a projection operator,
$P_\rho$, and a truncated Green function, $\widetilde G_\rho$, by
\begin{eqnarray}
\widehat{P_\rho f}(\xi)&= \chi_{T_\rho}(\xi)\cdot\widehat f(\xi)
\hbox{ and }\\ \label{eqn1.14}
(\widetilde G_\rho f)^{\wedge}(\xi)&=\chi_{T^C_\rho}(\xi)
\cdot[-\sigma(\xi)]^{-1} \widehat f(\xi) \label{eqn1.15}
\end{eqnarray}
for $f\in\cS(\R^n)$. Note that $\Delta_\rho\widetilde G_\rho=I-P_\rho$.

Choose a $\psi_3\in C^\infty_0(\R^{n-2})$, supported in $B^{n-2}(0;2)$,
radial and with $\psi_3\equiv 1$ on $\supp(\psi_1)$, and set
$\chi_3(x'')=\psi_3(\frac{x''-x''_0}{\delta})$. We now define the second
term, $u_1(x, \rho, \Pi)$ in the approximate solution by
\begin{equation}\label{eqn1.16}
u_1(x) = -\chi_3(x'')\widetilde G_\rho((\Delta_\rho+q(x))u_0(x))
\end{equation}
and set $u(x)=u_0(x)+u_1(x)$. Then $u_1$ (as well as $u_0$) is supported in
$\{x:\dist(x,\Pi)\le 2\delta\}$, yielding (\ref{eqn1.6}). We will see below
that $\|u_1\|_{L^2(\Omega)}\le C|\rho|^{-\epsilon}$ as $|\rho|\to\infty$, so
that
(\ref{eqn1.5}) holds as well, so
that the first part of (1.9) holds as well. To
start the proof of (\ref{eqn1.7}), note that
\begin{eqnarray*}
(\Delta_\rho+q)(u_0+u_1) & = & (\Delta_\rho+q)u_0-
(\Delta_\rho+q)\chi_3\widetilde G_\rho((\Delta_\rho+q)u_0) \\
& = & (\Delta_\rho+q)u_0-\chi_3(\Delta_\rho+q)\widetilde
G_\rho((\Delta_\rho+q)u_0) \\
&&\quad - [\Delta_\rho+ q,\chi_3]\widetilde G_\rho((\Delta_\rho+q)u_0) \\
& = & (\Delta_\rho+q)u_0 - \chi_3(I-P_\rho)(\Delta_\rho+q)u_0 -
\chi_3q\widetilde G_\rho(\Delta_\rho+q)u_0 \\
&&\quad - 2(\nabla\chi_3\cdot\nabla_{x''})\widetilde
G_\rho(\Delta_\rho+q)u_0-
(\Delta_{x''}\chi_3)\widetilde G_\rho(\Delta_\rho+q)u_0 \\
& = & \chi_3P_\rho(\Delta_\rho+q)u_0\\
&&\quad -[q\chi_3+2(\nabla \chi_3\cdot\nabla
_{x''})-(\Delta_{\chi''}\chi_3)]\widetilde G_\rho(\Delta_\rho + q)u_0
\end{eqnarray*}
on $\Omega$, since $\chi_3\equiv 1$ on $\supp(\chi_1)$. Now, since
$q_1\chi_3\in L^\infty$, $|\nabla\chi_3|\le C\delta^{-1}=c|\rho|^\beta$ and
$|\Delta_{x''}\chi_3|\le C\delta^{-2}=c|\rho|^{2\beta}$, (\ref{eqn1.7}) will
follow
if we can show that for some $\epsilon>0$,
\begin{eqnarray}
\|P_\rho(\Delta_\rho+q)u_0\|_{L^2(\Omega)} & \le & C|\rho|^{-\epsilon}, \\
\label{eqn1.17}
\|\,|D''|\widetilde G_\rho(\Delta_\rho+q)u_0\|_{L^2(\Omega)} & \le &
C|\rho|^{-\beta-\epsilon},\hbox{ and }\\ \label{eqn1.18}
\|\widetilde G_\rho(\Delta_\rho+q)u_0\|_{L^2(\Omega)} & \le &
C|\rho|^{-2\beta-\epsilon}, \label{eqn1.19}
\end{eqnarray}
with $C$ independent of $|\rho|>1$. Before proceeding to
prove these, we note that for any $u^{(1)}, u^{(2)}$ constructed in this way
for the same
two-plane $\Pi$,
\[
u^{(1)}_0(x)u^{(2)}_0(x)=\chi^2_0(x')\delta^{-(n-2)}\psi^2_1\left(
\frac{x''-x''_0}{\delta}\right)\to d\lambda_\Pi^\Omega\hbox{ in }\Omega
\]
as $\delta\to 0$ by  (1.11), while
$u^{(1)}_1u^{(2)}_0+u^{(1)}_0u^{(2)}_1 +  u^{(1)}_1u^{(2)}_1\to 0$ in
$L^2(\Omega)$, yielding (\ref{eqn1.8}). Thus, we are reduced to establishing
(1.17--1.19).

\section{$L^2$ estimates}

We will first prove (1.17)--(1.19) under the simplifying assumption that
$q_1,q_2\in C^{n-1+\sigma}(\overline{\Omega})$ for some $\sigma>0$, turning
to the Sobolev space case in Section 3.
Start by noting that the desired estimates (1.17)--(1.19) cannot be simply
obtained
from operator norms; for example, $\|P_\rho\|_{L^2\to L^2}=1$ for all
$\rho$. One needs to make use of the special structure of
$(\Delta_\rho+q)u_0$; we first deal with $\Delta_\rho u_0$, leaving
$q(x)\cdot u_0$ for the end.
So, we will show that $\|P_\rho\Delta_\rho u_0\|_{L^2}\le
C|\rho|^{-\epsilon}$, etc. Since $\nabla\chi_0\cdot\nabla\chi_1\equiv 0$,
\begin{equation}\label{eqn2.1}
\Delta_\rho
u_0=\chi_0\Delta_{x''}\chi_1+(\Delta_{x'}+4|\rho|
\overline\partial_{x'})(\chi_0)
\cdot\chi_1.
\end{equation}
The second term is supported on $\Omega^c$, but $P_\rho$ and $\widetilde
G_\rho$ are nonlocal operators and we need to control the contribution from
this term. However, because
$\Delta_{x'}(\chi_0)$ is a fixed, $\delta$-independent
element of $C^\infty_0(\R^2)$, this can be handled in the same way as the
$q(x)\cdot u_0$ terms of (1.17--1.19), which will be dealt with later.
The contribution from $4|\rho|\dbc\cdot\chi_1$ will be handled at the end.

\par So, for
the time being, we are interested in estimating
$\|P_\rho(\chi_0(x')\Delta_{x''}\chi_1(x''))\|_{L^2}$, etc. Now,
$\Delta_{x''}\chi_1(x'') = \delta^{-2}\chi_5(x'')$, where
$\chi_5(x'')=\delta^{-\frac{n-2}{2}}\psi_5
\left(\frac{x''-x''_0}{\delta}\right)$ is
associated with the radial function $\psi_5=\Delta_{x''}\psi_1$ as
$\chi_1$ is associated with $\psi_1$. Note for future use that
$\widehat\psi_5$ vanishes to second order at 0. Of course, $\chi_0\in
C^\infty_0\Rightarrow\widehat\chi_0\in\cS(\R^n)$, but looking ahead to
estimating the terms involving $q(x)\cdot u_0(x)$, we will now prove the
analogues of (1.17--1.19) where $P_\rho$ and $\widetilde G_\rho$ act on
$\chi_2(x')\Delta\chi_1(x'')$, under the weaker assumption that $\chi_2$ is
radial and satisfies the uniform decay estimate
\begin{gather*} \label{eqn2.2}
|\widehat\chi_2(\xi)|\le C(1+|\xi|)^{-\alpha} \tag*{$(2.2)_\alpha$}
\end{gather*}
for some $\alpha>0$.

\setcounter{equation}{2}
Now, by (\ref{eqn1.14}) and Plancherel,
\begin{eqnarray*}
\|P_\rho(\chi_2\Delta\chi_1)\|_{L^2(\Omega)} & \le &
\|(P_\rho(\chi_2\Delta\chi_1))^{\wedge}\|_{L^2(\R^n)} \\
& = &
\|\delta^{-2}|\widehat\chi_2(\xi')|\delta^{\frac{n-2}{2}}|\widehat\psi_5
(\delta\xi'')|\,\|_{L^2(T_\rho)}.
\end{eqnarray*}
The characteristic variety $\Sigma_\rho$, of which $T_\rho$ is a tubular
neighborhood, passes through the origin, and we may represent $\Sigma_\rho$
near $O$ as a graph over the $\xi''$-plane: $\Sigma_\rho=\Sigma^s_\rho\cup
\Sigma^n_\rho\cup\Sigma^e_\rho$, with
\begin{eqnarray}\label{eqn2.3}
\Sigma^s_\rho & = & \left\{\xi_1=0, \xi_2 =
|\rho|-(|\rho|^2-|\xi''|^2)^{\frac{1}{2}}, |\xi''|\le
\frac{|\rho|}{2}\right\}
\\
& \simeq & \left\{\xi_1=0, \xi_2=\frac{|\xi''|}{2|\rho|}, |\xi''|\le
\frac{|\rho|}{2}\right\}\nonumber
\end{eqnarray}
a neighborhood of the south pole $O$,
\begin{eqnarray}\label{eqn2.4}
\Sigma^n_\rho & = & \left\{\xi_1=0,
\xi_2=|\rho|+(|\rho|^2-|\xi''|^2)^{\frac{1}{2}}, |\xi''|\le \frac{|\rho|}{2}
\right\} \\
& \simeq & \left\{\xi_1=0, \xi_2=2|\rho|-\frac{|\xi''|^2}{2|\rho|},
|\xi''|\le \frac{|\rho|}{2}\right\} \nonumber
\end{eqnarray}
a neighborhood of the north pole $(0, 2|\rho|, 0, \ldots, 0)$,
and $\Sigma^e_\rho$ a neighborhood of the equator
$\{\xi\in\Sigma_\rho:\xi_2=|\rho|\}$. We have a corresponding decomposition
$T_\rho=T^s_\rho\cup T^n_\rho\cup T^e_\rho$, where, e.g.,
\begin{equation}\label{eqn2.5}
T^s_\rho\simeq\left\{(\xi', \xi''):\xi'\in B^2\left(\left(0,
\frac{|\xi''|^2}{2|\rho|}\right); |\rho|^{-\frac{1}{2}-\epsilon_0}\right),
|\xi''|\le
\frac{|\rho|}{2}\right\}.
\end{equation}
Recalling that $\chi_2$ and $\psi_3$ are radial, so are $\widehat\chi_2$ and
$\widehat\chi_3$, and by abuse of notation we consider these as functions of
one variable satisfying \ref{eqn2.2} and rapidly decreasing, respectively.
Thus, using polar coordinates in $\xi''$,
\begin{eqnarray}\label{eqn2.6}
\|\widehat{\chi_2\Delta\chi_1}\|^2_{L^2(T^s_\rho)} & \simeq &
\int^{\frac{|\rho|}{2}}_0 \int_{B^2\left(\left(0,
\frac{r^2}{2|\rho|}\right);|\rho|^{-\frac{1}{2}-\epsilon_0}\right)}
|\widehat\chi_2(\xi')|^2 d\xi' \delta^{n-6}|\widehat\psi_5
(\delta r)|^2r^{n-3}dr \nonumber \\
& \simeq & \int^{\sqrt{2}|\rho|^{\frac{1}{4}}}_0
\int_{B^2((0,0);|\rho|^{-\frac{1}{2}-\epsilon_0})}|\widehat\chi_2|^2d\xi'
\delta^{n-6}|
\widehat\psi_5(\delta r)|^2 r^{n-2}\frac{dr}{r} \\
&&+
\int^{\frac{|\rho|}{2}}_{\sqrt{2}|\rho|^{\frac{1}{4}}}\left|\widehat\chi_2
\left(\frac{r^2}{2|\rho|}\right)\right|^2 \cdot |B^2((0,0);
|\rho|^{-
\frac{1}{2}})|
\delta^{n-6}|\widehat\psi_5(\delta r)|^2r^{n-2}\frac{dr}{r}.\nonumber
\end{eqnarray}
Since we will be taking $\delta=|\rho|^{-\beta}$ with $\beta<\frac{1}{4}$,
if
we choose $0<\epsilon_0<2(\frac14-\beta)$, then
the quantity $|\rho|^{\frac{1}{4}}\delta\to\infty$ as $|\rho|\to\infty$ and
so
\begin{eqnarray}\label{eqn2.7}
\|\widehat{\chi_2\Delta\chi_1}\|^2_{L^2(T^s_\rho)} & \le &
c\frac{\delta^{-4}}{|\rho|^{1+2\epsilon_0}} \left(
\int^{\sqrt{2}|\rho|^{\frac{1}{4}}\delta}_{0}
|\widehat\psi_5(r)|^2r^{n-2} \frac{dr}{r} \right)\\
&&+\int_{\sqrt{2}|\rho|^{\frac{1}{4}}\delta}^{\frac{|\rho|}{2}\delta}
\left|\widehat\chi_2
\left(\frac{r^2}{2\delta^2|\rho|}\right)\right|^2
|\widehat\psi_5(r)|^2r^{n-2}\frac{dr}{r}
\nonumber \\
& \le & c(\delta^4|\rho|)^{-1},
\nonumber
\end{eqnarray}
which is $\le c|\rho|^{-2\epsilon}$ with $\epsilon=\frac12(1-4\beta)>0$.

The other contributions to $\|P_\rho\chi_2\Delta\chi_1\|_{L^2}$, coming from
$T^n_\rho$ and $T^e_\rho$ are handled similarly and are even smaller, due to
the decrease of $\widehat\chi_2$ and $\widehat\psi_5$.

We next turn to estimating $\|\,|D''|\widetilde G_\rho\Delta_\rho
u_0\|_{L^2}$; by the remark above, we may concentrate on the
$\chi_2\Delta\chi_1$ term of $\Delta_\rho u_0$. Then
\begin{equation}\label{eqn2.8}
\|\,|D''|\widetilde G_\rho(\chi_2\Delta\chi_1)\|^2_{L^2(\Omega)} \le
\|\,|\xi''|
(\sigma(\xi))^{-1}(\chi_2\Delta\chi_1)^{\wedge}(\xi)\|^2_{L^2(T^C_\rho)}.
\end{equation}

We may cover $T^C_\rho$ by $ T^{C,s}_\rho\cup T^{C,n}_\rho\cup
T^{C,e}_\rho\cup T^{C,\infty}_\rho$, where
\begin{equation}\label{eqn2.9}
T^{C,s}_\rho = \left\{\xi:\xi'\in B^2\left(\left(0,
\frac{|\xi''|^2}{2|\rho|}\right); |\rho|^{-\frac{1}{2}-\epsilon_0}\right)^C
\cap B^2
\left(\left(0,2|\rho| - \frac{|\xi''|^2}{2|\rho|}\right);
\frac{1}{4}|\rho|\right)^C, |\xi''|\le \frac{|\rho|}{2}\right\},
\end{equation}
$T^{C,n}_\rho$ is defined similarly,
\begin{equation}\label{eqn2.10}
T^{C,e}_\rho = \left\{\xi:\frac{|\rho|}{4}<\xi_2<\frac{7|\rho|}{4}, { }
|\rho|^{-\frac{1}{2}} < \dist(\xi,\Sigma_\rho)<|\rho|,
|\xi''|<2|\rho|\right\}
\end{equation}
and
\begin{equation}\label{eqn2.11}
T^{C,\infty}_\rho = \left\{\xi:|\xi|\ge 3|\rho|, |\xi''|\ge
\frac{3}{2}|\rho|\right\}.
\end{equation}
One has the lower bounds on $\sigma$,
\begin{equation}\label{eqn2.12}
|\sigma(\xi)|\ge \left\{ \begin{array}{ll}
C|\rho|\dist(\xi,\Sigma_\rho), & |\xi|\le 3|\rho| \\
C|\xi|^2, & |\xi|\ge 3|\rho|\end{array}\right.
\end{equation}
with $C$ (as always) uniform in $|\rho|$. The first inequality in
(\ref{eqn2.12}) follows from noting that
$\frac{1}{2}\nabla\sigma(\xi)=(\xi-|\rho|\vec{e_2})+i(|\rho|\vec{e_1})$, so
that $|\nabla\sigma(\xi)|=2\sqrt{2}|\rho|$ on ${\Sigma_\rho}$, while
the second follows from
$Re(\sigma(\xi))=\dist(\xi,|\rho|\vec{e_2})^2-|\rho|^2$. Using the first
estimate in (\ref{eqn2.12}), we can then dominate the contribution to the
right side of (\ref{eqn2.8}) from the region $T^{C,s}_\rho$ by
\begin{equation}\label{eqn2.13}
\delta^{n-6}\int_{|\xi''|\le\frac{|\rho|}{2}} \int_{B^2\left(\left(0,
\frac{|\xi''|^2}{2|\rho|}\right); |\rho|^{-\frac{1}{2}-\epsilon_0}
\right)^C} |\rho|^{-2}
\left|\xi' - \frac{|\xi''|^2}{2|\rho|}\vec{e_2}\right|^{-2}
|\widehat\chi_2(\xi')|^2d\xi' |\xi''|^2|\widehat\psi_5(\delta\xi'')|^2
d\xi''.
\end{equation}
The inner integral is the convolution
\[
\left. |\rho|^{-2}\left(|\widehat\chi_2|^2\ast_{\R^2}
\frac{\chi
\{|\xi'|\ge|\rho|^{-\frac{1}{2}-\epsilon_0}\}}{|\xi'|^2}\right)\right|_{\xi'
=
\frac{|\xi''|^2}{2|\rho|}\vec{e_2}}.
\]
An elementary calculation shows that, for $\widehat\chi_2$ satisfying
\ref{eqn2.2} for some $0<\alpha<1$, and any $0<a<1$,
\begin{equation}\label{eqn2.14}
|\widehat\chi_2|^2\ast_{\R^2} \frac{\chi\{|\xi'|\ge a\}}{|\xi'|^2} \le
\left\{ \begin{array}{ll}
C_1(1+\log(a^{-1})), & |\xi'|\le 1 \\
C_2|\xi'|^{-2}+C_3|\xi'|^{-2\alpha}\log\left(\frac{|\xi'|}{a}\right), &
|\xi'|\ge 1,\end{array}\right.
\end{equation}
so that, taking $a=|\rho|^{-\frac{1}{2}-\epsilon_0}$ and $|\xi'| =
\frac{|\xi''|^2}{2|\rho|}$,
the inner integral in (\ref{eqn2.13}) is
\[
\le \left\{
\begin{array}{ll}
C_1|\rho|^{-2}\log|\rho|, & 0<|\xi''|\le\sqrt{2}|\rho|^{\frac{1}{2}} \\
C_2|\xi''|^{-4} + C_3|\rho|^{2\alpha-2} |\xi''|^{-4\alpha} \log \left(
\frac{|\xi''|^2}{2|\rho|^{\frac{1}{2}-\epsilon_0}}\right), &
\sqrt{2}|\rho|^{\frac{1}{2}} \le |\xi''| \le \frac{|\rho|}{2}.
\end{array}
\right.
\]
Employing polar coordinates in $\xi''$ and rescaling by $\delta$, we see
that (\ref{eqn2.13}) is
\begin{eqnarray*}
& \le & C_1\delta^{-6}|\rho|^{-2}\log|\rho|
\int^{\sqrt{2}|\rho|^{\frac{1}{2}}\delta}_0
|\widehat\psi_5(r)|^2r^n\frac{dr}{r} \\
&& + C_2\delta^{-2}
\int^{\frac{|\rho|}{2}\delta}_{\sqrt{2}|\rho|^{\frac{1}{2}}\delta}
|\widehat\psi_5(r)|^2r^{n-4}\frac{dr}{r} \\
&& + C_3 \delta^{4\alpha-4}|\rho|^{2\alpha-2} \log|\rho|
\int^{\frac{|\rho|}{2}\delta}_{\sqrt{2}|\rho|^{\frac{1}{2}}\delta}
|\widehat\psi_5(r)|^2r^{n-2-4\alpha}\frac{dr}{r}.
\end{eqnarray*}
With $\delta=|\rho|^{-\beta}$, $\beta<\frac{1}{4}$,
$|\rho|^{\frac{1}{2}}\delta\to\infty$ as $|\rho|\to \infty$, and thus we
estimate this for any $N>0$ (using the rapid decay of $\widehat\psi_5$) by
\[
C_1|\rho|^{6\beta-2}\log|\rho|+C_2\delta^{-2}(|\rho|^{\frac{1}{2}}\delta)^{-
N}
+ C_3|\rho|^{(4-4\alpha)\beta+2\alpha-2} \log|\rho|(|\rho|^{\frac{1}{2}}
\delta)^{-N},
\]
the first term of which will be less than the desired
$|\rho|^{-2\beta-2\epsilon}$,
for any $\alpha>0$, if $\beta<\frac{1}{4}$ and $\epsilon=
\frac12(1-4\beta)$; the second and third
terms are rapidly decaying simply because $\beta<\frac{1}{2}$.

Moving ahead for the moment to (\ref{eqn1.19}), the contribution to
$\|\widetilde G_\rho \chi_2\Delta\chi_1\|^2_{L^2}$ (which we want $\le
C|\rho|^{-4\beta-2\epsilon}$) from $T^{C,s}_\rho$ is handled in the same
fashion, the only differences being the absence of the multiplier
$|D''|^\wedge = |\xi''|$ on the left and the improved gain we are demanding
on the right. Taking these into account, we need to control
\begin{multline}
C_1\delta^{-4}|\rho|^{-2}\log|\rho|\int^{\sqrt{2}|\rho|^{\frac{1}{2}}\delta}
_0
|\widehat\psi_5(r)|^2r^{n-2}\frac{dr}{r} \\
+ C_2 \int^{\frac{1}{2}|\rho|\delta}_{\sqrt{2}|\rho|^{\frac{1}{2}}\delta}
|\widehat\psi_5(r)|^2r^{n-6}\frac{dr}{r} \\
+ C_3\delta^{4\alpha-2}|\rho|^{2\alpha-2}\log|\rho|
\int^{\frac{1}{2}|\rho|\delta}_{\sqrt{2}|\rho|^{\frac{1}{2}}\delta}
|\widehat\psi_5(r)|^2 r^{n-4-4\alpha} \frac{dr}{r} \\
\le C_1\delta^{-4}|\rho|^{-2}\log|\rho| +
C_2(|\rho|^{\frac{1}{2}}\delta)^{-N}
+ C_N\delta^{4\alpha-2} |\rho|^{2\alpha-2}
\log|\rho|(|\rho|^{\frac{1}{2}}\delta)^{-N},
\end{multline}
and this is $\le C|\rho|^{-4\beta-2\epsilon}$ provided $\beta<\frac{1}{4}$,
$\epsilon<
\frac12(1-4\beta)$
and $N$ is sufficiently large.

The contributions to (1.18)  from $T^{C,n}_\rho$
and $T^{C,e}_\rho$ are handled similarly. To treat the contribution from
$T^{C,\infty}_\rho$, we use the second estimate in (\ref{eqn2.12}) and
calculate (for (1.18)
\begin{multline}
\|\,|\xi''|(\sigma(\xi))^{-1}(\chi_2\Delta\chi_1)^\wedge(\xi)\|^2_{L^2(T^{C,
\infty}_\rho)} \\
\le C\int\kern-5pt\int_{|\xi|\ge 3|\rho|} \delta^{n-6}
|\widehat\chi_2(\xi')|^2|\widehat\psi_5(\delta\xi'')|^2
\frac{|\xi''|^2d\xi'd\xi''}{|\xi|^4} \\
\le C\left(\int_{|\xi''|\le|\rho|}\delta^{n-6}|\rho|^{-2\alpha-2}
|\widehat\psi_5(\delta\xi'')|^2|\xi''|^2d\xi''\right. \\
\left.+ \int_{|\xi''|\ge|\rho|} \delta^{n-6}
|\widehat\psi_5(\delta\chi'')|^2
|\xi''|^{-2\alpha} d\xi''\right) \\
= C\left(\delta^{-6}|\rho|^{-2\alpha-2} \int^{|\rho|\delta}_{0}
|\widehat\psi_5(r)|^2r^n\frac{dr}{r}\right. \\
\left.+ \delta^{2\alpha-4}\int^\infty_{|\rho|\delta}
|\widehat\psi_5(r)|^2r^{n-2-2\alpha}\frac{dr}{r}\right) \\
\le C(\delta^{-6}|\rho|^{-2\alpha-2} +
\delta^{2\alpha-4}(|\rho|\delta)^{-N}),\qquad \forall\;N>0,
\end{multline}
which, for $\delta=|\rho|^{-\beta}$ and $N$ large is $\le
C|\rho|^{-2\beta-2\epsilon}$ provided $\beta<\frac{1}{4}$
and $\epsilon<\alpha+1-4\beta$.
A similar analysis holds for the $T^{C,\infty}_\rho$ contribution to
(1.19).

We now turn to controlling the $q(x)u_0(x)$ terms in
(1.17)--(1.19), as well as the contributions from the
$\Delta(\chi_0)\cdot\chi_1$ term in (\ref{eqn2.1}). Note that since
$q(x)$ is $C^{n-1+\sigma}$ (for some $\sigma>0$), $q(x)$ has an
extension (see, e.g., [St70,Ch.6]) to a $C^{n-1+\sigma}$ function of compact
support on
$\R^n$, which we also denote by $q$. The restriction of $q$ to any $\Pi\in
M_{2,n}$ is still
$C^{n-1+\sigma}$.

Let $\{D_t:0<t<\infty\}$ be the one-parameter group of partial dilations on
$\cS'(\R^{n^\ast})$,
\[
(D_tf)(\xi', \xi'')=t^{n-2}f(\xi', t\xi''),
\]
which, for $f, g\in L^1$, satisfy $\int_{\R^n}D_tfd\xi=\int_{\R^n}fd\xi$ and
$D_t(f\ast g)=D_tf\ast D_tg$. Then
\begin{eqnarray}\label{eqn2.15}
\widehat{qu}_0(\xi) & = & \widehat q\ast\widehat
u_0(\xi)=D_\delta(D_{\delta^{-1}}
\widehat q)\ast\delta^{-\frac{n-2}{2}}D_\delta(\widehat\chi_0(\xi')
\widehat\psi_1(\xi'')e^{ix''_0\cdot\xi''}) \\
& = & D_\delta(D_{\delta^{-1}}(\widehat q) \ast
\delta^{-\frac{n-2}{2}}\widehat\chi_0 \widehat\psi_1 e^{ix''_0\cdot\xi''}).
\nonumber
\end{eqnarray}
Now, as $\delta=|\rho|^{-\beta}\to 0$, $D_{\delta^{-1}}(\widehat q) =
\delta^{-(n-2)} \widehat q(\xi', \delta\xi'')$ converges weakly
 to the singular measure
\begin{equation}\label{eqn2.16}
Q(\xi')\otimes\delta(\xi'') = Q(\xi')d\xi',
\end{equation}
where $Q(\xi') = \int_{\R^{n-2}}\widehat q(\xi', \xi'')d\xi''$;
note that $q\in C^{n-1+\gamma}$ implies that the integral defining $Q$
converges and
$Q$ satisfies  (2.2)$_{1+\gamma}$. Letting
$F(\xi)=\widehat\chi_0(\xi')\widehat\psi_1(\xi'') e^{ix''_0\cdot\xi''}$, it
follows from (\ref{eqn2.15}) that
\begin{eqnarray}\label{eqn2.17}
\hbox{  }\widehat{qu}_0(\xi) & = & D_\delta(D_{\delta^{-1}}(\widehat q) \ast
\delta^{-\frac{n-2}{2}}F) \\
& = & D_\delta((Qd\xi')\ast\delta^{-\frac{n-2}{2}}F) +
D_\delta((D_{\delta^{-1}}
\widehat q-Qd\xi')\ast\delta^{-\frac{n-2}{2}}F).\nonumber
\end{eqnarray}
If we define $\widehat\chi_4(\xi')=Q\ast_{\R^2}\widehat\chi_0(\xi')$, then
$\widehat\chi_4$ also satisfies condition (2.2)$_{1+\gamma}$ (and thus
(2.2)$_{\alpha'}$
for $0<\alpha'<1$, so that (2.14) can be applied), and the first term
in (\ref{eqn2.17}) is
\begin{equation}\label{eqn2.18}
D_\delta((Qd\xi')\ast\delta^{-\frac{n-2}{2}}F) = \widehat\chi_4(\xi')
\delta^{\frac{n-2}{2}}\widehat\psi_1(\delta\xi'') e^{i\delta
x''_0\cdot\xi''}.
\end{equation}

Thus, the contributions to $\|P_\delta(qu_0)\|_{L^2}$, $\|\,|D''|\widetilde
G_\rho(qu_0)\|_{L^2}$ and $\|\widetilde G_\rho(qu_0)\|_{L^2}$  from the
first term in (\ref{eqn2.17}) may be handled as the main
$\chi_2\Delta\chi_1$ term was earlier, with the obvious absence of the
factor $\delta^{-2}$. To control the contributions from the second term in
(\ref{eqn2.17}), we use the elementary

\begin{lem}
Let $\varphi(x)$, $f(x)$ be functions on $\R^m$ such that $\varphi(x)$,
$|x|\varphi(x)$, $f(x)$ and $|\nabla f(x)|$ are in $L^1(\R^m)$. Then,
$\forall\;\epsilon>0$
\begin{eqnarray}\notag
\left|\left(\epsilon^{-m}\varphi\left(\frac{x}{\epsilon}\right) \right.
\right.& - &
\left.\left. \left(\int_{\R^m}\varphi dy\right)\delta(x)\right)\ast
f(x)\right|\\
&\le & C_m(\|\varphi\|_{L^1}+\|\,|x|\varphi\|_{L^1})\cdot
(\|f\|_{L^\infty(B(0;|x|-1))} +
\|\nabla f\|_{L^\infty (B(x;1))})\cdot\epsilon.\nonumber
\end{eqnarray}
\end{lem}

Applying this for $\epsilon=\delta$, $\xi'\in\R^2$ fixed, and using
$F\in\cS$, $|\widehat q(\xi)|\le C(1+|\xi|)^{-(n-1+\gamma)}$, we find that,
$\forall\;N>0$
\begin{equation}\label{eqn2.20}
|(D_{\delta-1}(\widehat q) - Qd\xi') \ast F(\xi)|\le C_N(1+|\xi'|)^{-\gamma}
(1+|\xi''|)^{-N}\delta.
\end{equation}
Hence, the second term in (\ref{eqn2.17}) is $\le
C_N\delta^{\frac{n}{2}}(1+|\xi'|)^{-\gamma}(1+|\delta\xi''|)^{-N}$ and this
allows the contributions to (1.17)--(1.19) to be dealt
with as the $\chi_2\Delta_{x''}\chi_1$ term was before.
\par Finally, we need to establish the estimates (1.17--1.19) for the
$4|\rho|\dbc$ term in (\ref{eqn2.1}); thus, we need to show
\begin{eqnarray}
\|P_\rho\bigl(\dbc\cdot\chi_1\bigr)\|_{L^2} & \le & C|\rho|^{-1-\epsilon},
\\
\label{eqn2.23}
\|\,|D''|\widetilde G_\rho\bigl(\dbc\cdot\chi_1\bigr)\|_{L^2} & \le &
C|\rho|^{-1-\beta-\epsilon},\hbox{ and }\\ \label{eqn2.24}
\|\widetilde G_\rho\bigl(\dbc\cdot\chi_1\bigr)\|_{L^2} & \le &
C|\rho|^{-1-2\beta-\epsilon}, \label{eqn2.25}
\end{eqnarray}
for some $\epsilon>0$. Using the fact that $\widehat\dbc(\xi')$ is
rapidly decreasing and vanishes to first order at $\xi'=0$, we may
replace (\ref{eqn2.6}) with

\begin{eqnarray}\nonumber
\|\widehat{\dbc\chi_1}\|^2_{L^2(T^s_\rho)} & \simeq &
\int^{\frac{|\rho|}{2}}_0 \int_{B^2\left(\left(0,
\frac{r^2}{2|\rho|}\right);|\rho|^{-\frac{1}{2}-\epsilon_0}\right)}
|\widehat\dbc(\xi')|^2 d\xi' \delta^{n-2}|\widehat\psi_1
(\delta r)|^2r^{n-3}dr \nonumber \\
& \le c_N\Bigl( & \int^{\sqrt{2}|\rho|^{\frac{1-2\epsilon_0}{4}}}_0
|\rho|^{-2-4\epsilon_0}\delta^{n-2}|
\widehat\psi_1(\delta r)|^2 r^{n-2}\frac{dr}{r} \nonumber\\
&&+
\int^{\sqrt{2}|\rho|^{\frac12}}
_{\sqrt{2}|\rho|^{\frac{1-2\epsilon_0}{4}}}
(\frac{r^2}{2|\rho|})^2|\rho|^{-1-2\epsilon_0}
\delta^{n-2}|\widehat\psi_1(\delta r)|^2r^{n-2}\frac{dr}{r}\nonumber\\
&&+
\int_{\sqrt{2}|\rho|^{\frac12}}
^{\frac{|\rho|}2}
(\frac{r^2}{2|\rho|})^{-N}|\rho|^{-1-2\epsilon_0}
\delta^{n-2}|\widehat\psi_1(\delta r)|^2r^{n-2}\frac{dr}{r}\Bigr)
\nonumber \\
&\le c_N\Bigl(&
|\rho|^{-2-4\epsilon_0}
\int^{\sqrt{2}|\rho|^{\frac{1-2\epsilon_0}{4}}\delta}_0
|\widehat\psi_1|^2 r^{n-2}\frac{dr}r\nonumber\\
&&+
|\rho|^{-3-2\epsilon_0}\delta^{-4}
\int^{\sqrt{2}|\rho|^{\frac12}}
_{\sqrt{2}|\rho|^{\frac{1-2\epsilon_0}{4}}\delta}
|\widehat\psi_1|^2 r^{n+2}\frac{dr}r\nonumber\\
&&+
|\rho|^{-1-2\epsilon_0+N}\delta^{2N}
\int_{\sqrt{2}|\rho|^{\frac12}\delta}
^{\frac{|\rho|}2\delta}
|\widehat\psi_1|^2 r^{n-2-2N}\frac{dr}r\Bigr)\nonumber\\
&\le c_N &\bigl( |\rho|^{-2-4\epsilon_0}
+|\rho|^{-3-2\epsilon_0+4\beta}
(|\rho|^{\frac{1-2\epsilon_0}4-\beta})^{-N'}\nonumber \\
&&+|\rho|^{-1-2\epsilon_0+N-2N\beta-N'(\frac12-\beta)}\bigr)\nonumber\\
\end{eqnarray}
for any $N,N'\ge 0$.  As before, the contributions
from $T^n_\rho$ and $T^e_\rho$ are handled similarly.
Since $\epsilon_0<\frac12-2\beta$, if $N'$ is chosen large enough
 this yields (2.23) with
$\epsilon\le 2\epsilon_0$, which is weaker than the previously
imposed $\epsilon<\frac12(1-4\beta)$.
\par The desired estimates (2.23),(2.24) are even easier and hold for
any $\beta<\frac12$. The contribution to (2.24) from $T^{C,s}_\rho$ is
controlled as in (2.13), but with the factor $\delta^{n-2}$
and with the $\widehat{\chi_2}$ in the integrand replaced by $\widehat\dbc$;
this
is then dominated in the same manner as below (2.14). The $T^{C,s}_\rho$
contribution to (2.25) is estimated as in (2.15), but with
the absence of the $\delta^{-4}$.
All other contributions are dealt with similarly.

This concludes the proof of Thm.\ref{thm4} for the case of potentials
in the H\"older class $C^{n-1+\sigma}(\overline{\Omega}), \sigma>0$.
The restrictions on $\beta$ and $\epsilon$ that we ahve needed are
that $\beta<\frac14$ and $\epsilon<\frac12(1-4\beta)$.

\section{Remarks}

(i) The proof of Thm. 4 needs to be
slightly modified if we assume that the potential
$q(x)$ belongs to the Sobolev space $H^{\frac{n}2+\sigma}({\Omega})$ for
some $\sigma>0$. Since $\po$ is Lipschitz, such a $q(x)$ can, by the
Calder\'on extension
theorem, be extended to be in $H^{\frac{n}2+\sigma}(\R^n)$. Again denoting
the extension
by $q$, one has by Cauchy-Schwarz

\begin{equation}\label{eqn3.1}
\int_{\R^2}\Bigl( \int_{\R^{n-2}} (1+|\xi''|)
|\hat q(\xi',\xi'')| d\xi''\Bigr)^2
(1+|\xi'|)^\sigma d\xi' \le c
\bigl(\|q\|_{\frac{n}2+\sigma}\bigr)^2
\end{equation}

Thus, $Q$ as in (2.18) belongs to $L^2(\R^2;(1+|\xi'|)^\sigma d\xi')$, so
that $\widehat\chi_4= Q*_{\R^2}\widehat\chi_0\in L^2(\R^2;(1+|\xi'|)^\sigma
d\xi')
\cap L^\infty$. Replacing the uniform decay estimate (2.2)$_\alpha$ with

\begin{gather*} \label{eqn3.2}
\widehat\chi_2\in L^2(\R^2;(1+|\xi'|)^\sigma d\xi')  \tag*{$(3.2)_\sigma$}
\end{gather*}
\setcounter{equation}{2}
will allow us to handle the first term in (2.19). Furthermore,
if for $\xi'$ fixed, we let $\phi(\cdot)=\widehat{q}(\xi',\cdot)$ in
Lemma 5, then $\phi(\xi'')$ and $|\xi''|\phi(\xi'')$ are in
$L^1(\R^{n-2})$ with norms (as functions of $\xi'$) in
$L^2(\R^2;(1+|\xi'|)^\sigma d\xi')$, and so the second term in (2.19)
is $\le c_N\widehat{\chi_6}(\xi')(1+|\delta\xi''|)^{-N}, \forall N$, with
$\widehat{\chi_6}$ satisfying condition (3.2)$_\sigma$. So, we are reduced
to repeating the analysis of Section 2 with (2.2)$_\alpha$ replaced by
(3.2)$_\sigma$.
The decay of $\widehat{\chi_2}$ was used in only two places in the
argument. In (2.14), under (3.2)$_\sigma$, we have the same estimate except
for
the absence of $|\xi'|^{-2\alpha}$; however, this loss is absorbed into
terms
rapidly decreasing in $|\rho|^{\frac12}\delta=|\rho|^{\frac12-\beta}$ where
(2.14) is used. On the other hand, in (2.16) we may estimate the inner
integral
by

\begin{eqnarray}\label{eqn3.3}
\int_{|\xi'|\ge
2|\rho|}|\widehat{\chi_2}(\xi')|^2\frac{d\xi'}{(|\xi'|^2+|\xi''|^2)^2}
& \le \int_{\R^2}|\widehat{\chi_2}|^2\frac{d\xi'}{(1+|\xi'|)^\sigma
|\xi'|^4}\\
& \le c|\rho|^{-4-\sigma} \hbox{ if } |\xi''|\le\rho \nonumber
\end{eqnarray}
and
\begin{equation}\label{eqn3.4}
\int_{\R^2}|\widehat{\chi_2}(\xi')|^2\frac{d\xi'}{(|\xi'|^2+|\xi''|^2)^2}
\le c|\xi''|^{-4}\hbox{ if } |\xi'|\ge\rho,
\end{equation}
so that

\begin{eqnarray}\label{eqn3.5}
&\|\,|\xi''|
(\sigma(\xi))^{-1}(\chi_2\Delta\chi_1)^\wedge(\xi)\|^2_{L^2(T^{C,\infty}_
\rho)} \\
&
\le C\left(\int_{|\xi''|\le|\rho|}\delta^{n-6}|\rho|^{-4-\sigma}
|\widehat\psi_5(\delta\xi'')|^2|\xi''|^2d\xi''
 + \int_{|\xi''|\ge|\rho|} \delta^{n-6} |\widehat\psi_5(\delta\chi'')|^2
|\xi''|^{-2} d\xi''\right) \nonumber \\
&= C\left(\delta^{-6}|\rho|^{-4-\sigma} \int^{|\rho|\delta}_{0}
|\widehat\psi_5(r)|^2r^n\frac{dr}{r}\right. \nonumber \\
&\left.+ \delta^{-2}\int^\infty_{|\rho|\delta}
|\widehat\psi_5(r)|^2r^{n-4}\frac{dr}{r}\right) \nonumber \\
& \le C_N(\delta^{-6}|\rho|^{-4-\sigma} + \delta^{-2}(|\rho|\delta)^{-N})
\nonumber \\
& = C_N\left( |\rho|^{6\beta-4-\sigma}+
|\rho|^{2\beta}(|\rho|^{\beta-\frac12})^N\right) ,\qquad \forall N,
\nonumber
\end{eqnarray}
which is $\le c|\rho|^{-2\beta-\epsilon}$ for $N$ sufficiently large, since
$\beta<\frac12$. The restrictions on $\beta$ and $\epsilon$ are as before.

(ii) The construction of the approximate solutions given by Thm. 4 may be
generalized
by taking $\chi_0$ to be an arbitrary analytic function of $z=x_1+ix_2$,
defined on a domain
$\Pi\cap\Omega\subset\subset\Omega'\subset\Pi$. Since
$\dbc=\Delta_{x'}\chi_0\equiv 0$
on $\Omega$, the resulting $u= u_0+u_1$ is still an approximate solution in
the sense
of Thm. 4, except that (1.8) no longer applies. Thus, Thm. 1 can be
strengthened
to conclude that $(q_1-q_2)|_\Pi$ is orthogonal in
$L^2(\Pi\cap\Omega,d\lambda_\Pi)$ to
the Bergman space $A^2(\Pi\cap\Omega)$ of square-integrable holomorphic
functions
on $\Pi\cap\Omega$. Furthermore, by repeating the construction using
$\overline{\rho}=
\frac1{\sqrt{2}}|\rho|(\omega_R-i\omega_I)$, which induces the conjugate
complex structure
on $\Pi$, for which the $\overline{\partial}$ operator equals the $\partial$
operator induced by $\rho$, we obtain
that $(q_1-q_2)|_\Pi$ is also orthogonal to the conjugate Bergman
space $\overline{A}^2(\Pi\cap\Omega)$ of anti-holomorphic functions. (The
analogue of this
in two dimensions was obtained in \cite{SU87b}.)
It would be interesting to make further use of this information.

(iii) To obtain variants of Thm. 1 establishing smaller sets of uniqueness
in $\po$, it might be
useful to use approximate solutions associated to different two-planes. For
this, it
seems necessary to construct approximate solutions with much thinner
supports, i.e., to
overcome the restriction $\beta<\frac14$ in Thm. 4. Such an improvement
 might also be useful in
extending the results to $q_j\in L^\infty$.

\end{document}